\documentstyle[11pt]{article}

\def\be{\begin{equation}}
\def\ee{\end{equation}}
\def\bea{\begin{eqnarray}}
\def\eea{\end{eqnarray}}

\def\1{\'{\i}}

\def\R{{\rm I\kern-.2em R}}

\def\oper{{\cal O}}

\def\ados{\lambda}

\def\luisH{{\cal H}}
\def\luisP{{\cal P}}
\def\luisM{{\cal M}}
\def\luisK{{\cal K}}
\def\luisD{{\cal D}}
\def\luisC{{\cal C}}

\begin{document}

\baselineskip=18pt

 \
\hfill \today

\
\vspace{1.5cm}

 \begin{center}

{\LARGE{\bf{On quantum algebra symmetries}}}

\smallskip

{\LARGE{\bf{of discrete Schr\"odinger equations}}}

\end{center}

\bigskip\bigskip

\begin{center} 
A. Ballesteros$\dagger$, F.J. Herranz$\dagger$, J.
Negro$\ddagger$ and L.M. Nieto$\ddagger$ 
\end{center}

\begin{center} 
{\it ${\dagger}$ Departamento de F\1sica\\
Universidad de Burgos,   E-09001  Burgos, Spain} 
\end{center}

\begin{center} 
{\it ${\ddagger}$ Departamento de F\1sica
Te\'orica\\
Universidad de Valladolid, E-47011  Valladolid, Spain}
\end{center}

\bigskip\bigskip\bigskip

\begin{abstract}\baselineskip=18pt
Two non-standard quantum deformations of the (1+1) 
Schr\"odinger algebra are identified with the symmetry
algebras of  either a space or time uniform lattice
discretization of the Schr\"odinger equation. For both
cases, the deformation parameter of the corresponding
Hopf algebra can be interpreted as the step of  the 
lattice. In this context, the introduction of nonlinear 
maps defining Schr\"odinger and $sl(2,\R)$ quantum algebras
with classical commutation rules turns out to be relevant. 
The problem of finding a quantum algebra linked to the full
space-time discretization is also discussed.
\end{abstract}

\bigskip\bigskip



%%%%%%%%%%%%%%%%%%%%%%%%%%%%%%%%%%%%%%%%%%%%%%

\newpage

\section*{I. INTRODUCTION}

Since they were introduced, quantum algebras have
been connected with different versions of space-time
lattices through several algebraic
constructions that have no direct relationship with
the usual Lie symmetry theory.${}^{1-3}$ Recent 
works${}^{4,5}$ have also developed new techniques for
dealing with the symmetries of difference or
differential-difference equations and have tried to adapt
in this field the standard methods that have been so
successful with differential equations. An exhaustive study
for the discretization on $q$-lattices of classical linear
differential equations has shown that their symmetries
obeyed to
$q$-deformed commutation relations with respect to the Lie
algebra structure of the continuous symmetries.${}^{6-8}$
However, Hopf
algebra structures underlying these $q$-symmetry  algebras 
have been not found.

When the discretization of linear
equations is made on uniform lattices it is well known that
the relevant symmetries preserve the Lie algebra
structure.${}^{9,10}$
Perhaps, this is the reason why 
the symmetry approach to these equations has  never been
related to quantum algebras. We will address this question 
in this paper for the discrete $(1{+}1)$-Schr\"odinger
equation. In this case we will show that Lie algebras of
(discrete) symmetries  can be put in correspondence to
non-trivial quantum Hopf algebra structures of the
non-standard (or triangular) type.${}^{11-16}$

To begin with, let us consider the following discrete
version  of the  heat or (time imaginary) Schr\"odinger
equation (SE) on a two-dimensional uniform lattice${}^{9}$
\be
(\Delta_x^2- 2 m \Delta_t)\phi(x,t)=0 .
\label{za}
\ee
The  difference operators $\Delta_x$ and $\Delta_t$
which appear in (\ref{za}) can be expressed in terms of 
shift operators
$T_x= e^{\sigma \partial_x}$ and $T_t=e^{\tau \partial_t}$ 
as
\be
\Delta_x=\frac{T_x-1}{\sigma}%\label{zb}
\qquad
\Delta_t=\frac{T_t-1}{\tau} \label{zc},
\ee
where the  parameters
$\sigma$ and $\tau$  are the lattice constants  in the 
space
$x$ and time $t$ directions, respectively.  The action of
$\Delta_x$ or
$\Delta_t$ on a function
$\phi(x,t)$ consists in a discrete derivative, which in 
the limit
$\sigma\to 0$ and $\tau\to 0$ come into $\partial_x$ and
$\partial_t$, respectively.

We will say that a given operator $\oper$ is a symmetry of
the linear equation
$E
\phi(x,t)=0$ if $\oper$ transforms solutions into  
solutions, that is, if $\oper$ is such that
\be
 E\, \oper = \Lambda\, E, \label{zf}
\ee
 where $\Lambda$ is
another operator. In this way,  the symmetries of the 
equation (\ref{za})  were computed,${}^{9}$ showing that
they realized the Schr\"odinger algebra ${\cal S}$, which 
is exactly the same result as for the
continuous case.${}^{17,18}$

On the other hand, two quantum deformations of the 
Schr\"odinger algebra have been recently
obtained,${}^{19,20}$ both of them endowed with a 
triangular Hopf structure.  The former was derived by
starting from the non-standard quantum deformation of its
harmonic oscillator subalgebra
$h_4$, while the latter was constructed by means of the 
non-standard quantum deformation of the extended
${\overline{sl}}(2,\R)$ subalgebra.  In Sec.~II  we
show that the symmetry algebra${}^{9}$ of
the space discretization of the SE obtained from (\ref{za})
by taking the limit $\tau\to 0$ is just
the quantum Schr\"odinger algebra given by
Ballesteros {\it et al.\/},${}^{19}$ and  the deformation
parameter
$z$ is related with the  space lattice constant $\sigma$.  
Likewise, we also show in Sec.~III  that the time
discretization of the SE obtained from (\ref{za}) by means
of the limit $\sigma\to 0$ has the quantum Schr\"odinger
algebra obtained by Ballesteros {\it et al.\/}${}^{20}$ as
its symmetry algebra; in this case, the time lattice step
$\tau$ plays the role of the deformation parameter. As a
consequence, a relationship between non-standard
deformations and regular lattice discretizations can be
established. In order to derive both results, the
introduction of nonlinear maps transforming the
aforementioned quantum algebras into non-cocommutative Hopf
algebras with non-deformed commutation rules turns out to
be essential. Such nonlinear transformations are  
explicitly given and, in the case of the quantum
algebra just mentioned,${}^{20}$ it is used to
derive a new nonlinear map for the jordanian quantum
deformation of
${{sl}}(2,\R)$. Finally, the problem of finding a quantum
algebra related to a full space-time discretization of the
SE is   discussed in Sec.~IV.

%%%%%%%%%%%%%%%%%%%%%%%%%%%%%%%%%%%%%%%%%%%%%%

\section*{II. A DISCRETE SPACE SCHR\"ODINGER EQUATION}

We will use the familiar notation for  the Schr\"odinger
generators: time translation $H$, space translation $P$,
Galilean boost $K$, dilation
$D$, conformal transformation $C$, and the central
generator $M$.${}^{17,18}$
Let $U_z^{(s)}({\cal S})$ be the quantum
Schr\"odinger algebra obtained by Ballesteros {\it et
al.\/}${}^{19}$ whose underlying Lie bialgebra is generated
by the non-standard classical $r$-matrix
\be
r^{(s)}=z
P\wedge (D +\frac 12 M).
\label{baa}
\ee
The coproduct of
$U_z^{(s)}({\cal S})$ is
\be
\begin{array}{l}
\Delta(M)=1\otimes M + M
\otimes 1 \\  [2pt]
\Delta(P)=1\otimes P + P \otimes 1  \\  [2pt]
\Delta(H)=1\otimes H + H\otimes e^{-2zP} \\  [2pt]
\Delta(K)=1\otimes K + 
K\otimes e^{zP}-z D'\otimes e^{zP}M \\  [4pt]
\Delta(D)=1\otimes D + D\otimes e^{zP} + \frac 12
M\otimes(e^{zP}-1) \\  [4pt]
\Delta(C)=1\otimes C + C\otimes  e^{2zP}+
\frac z2 K\otimes  e^{zP} D'
-\frac z2  D'\otimes  e^{zP}(K+z D' M) , \end{array}  
\label{ba}
\ee
and its commutation rules are
\be
\begin{array}{l}
[D,P]=\frac 1z ({1- e^{z P}}) \qquad\    [D,K]=K \qquad\
[K,P]=M  e^{z P }\qquad\  [M,\,\cdot\,]=0 \\  [4pt]
[D,H]= - 2 H \qquad\qquad   [D,C]=2 C -\frac z2 K D' 
\qquad\qquad [H,P]=0 \\  [4pt]
[H,C]= \frac 12(1+  e^{-z P}) D' -\frac 12 M +  z K H
\qquad\qquad  [K,C]=-\frac z2 K^2  \\  [4pt]
[P,C]=-\frac 12 (1+ e^{zP}) K - \frac z2   e^{zP} M D'
\qquad\qquad   [K,H]=\frac 1z({1-e^{-zP}})
\end{array}
\label{bb}
\ee
where hereafter we shall use the notation  
$D'=D+\frac 12 M$
in order to simplify some expressions.  The following
differential-difference realization of (\ref{bb}) can be 
found:
\bea
&& P=\partial_x \qquad\qquad H=\partial_t \qquad  \qquad M=
m \cr
&& K=- t \,\frac {{1- e^{-z \partial_x}}}{z}   - m x e^{z
\partial_x}
\qquad D=2  t  \partial_t 
 + x\, \frac {{  e^{z \partial_x}-1}}{z} +
\frac 12  \cr
&& C= t^2 \partial_t  e^{-z \partial_x} + t x \left(
\frac  {\sinh z\partial_x}{z} 
- z m \partial_t  e^{z \partial_x}\right)
+\frac 12 m x^2 e^{z \partial_x}\cr
&& \qquad -\frac 14 t \left\{  1- 3
e^{-z \partial_x} +
  m (1-e^{-z \partial_x})\right\}  
 +\frac 14 z m (1-   m) x  e^{z
\partial_x} .
\label{bc}
\eea

Note that the Galilei
generators
$\{K,H,P,M\}$ close a deformed subalgebra (but not a Hopf
subalgebra) whose Casimir  is
\be
E_z^{(s)} = \left(\frac {1-e^{- z P}}{z} \right)^2 - 2 M H.
\label{cass}
\ee 
The action of
$E_z^{(s)}$ on $\phi(x,t)$ through (\ref{bc}) provides a  
space discretization of the SE:
\be
\left(     \left(\frac {1-e^{-z \partial_x} } {z}\right)^2
- 2 m \partial_t\right) \phi(x,t)=0 .
\label{bd}
\ee
Furthermore, the
quantum Hopf algebra $U_z^{(s)}({\cal S})$ is a symmetry 
algebra of this equation since
\bea
&&[E_z^{(s)},X]=0\quad \mbox{for}\quad
X\in\{K,H,P,M\}\qquad
[E_z^{(s)},D]=  2 E_z^{(s)} \cr
&&[E_z^{(s)},C]=\left\{ t (e^{-z \partial_x} +1) - z m x 
e^{z
\partial_x} \right\} E_z^{(s)} .
\label{be}
\eea

Now we will relate these results with  the discretizations 
of the SE studied by Floreanini {\it et
al.\/}${}^{9}$  through the non-linear change of basis
defined by
\bea
&& \luisH=H\qquad \luisP=\frac   {1-e^{-zP}}{z}
\qquad \luisM=M\cr
&& \luisD=D+\frac 12 (1-e^{zP})\qquad \luisK=-2K
-z M e^{zP}\cr
&& \luisC=C-\frac z2 K D' + \frac z2 K
e^{zP}-\frac{z^2}8 M e^{2 zP}  \qquad \sigma=-z .
\label{bf}
\eea
 In this new basis the commutation
rules (\ref{bb}) of the Hopf algebra  $U_\sigma^{(s)}({\cal
S})$ are
\be
\begin{array}{llll}
 [\luisD,\luisP]=-\luisP \quad  &[\luisD,\luisK]=\luisK 
\quad &[\luisK,\luisP]=-2\luisM \quad
&[\luisM,\,\cdot\,]=0\cr
 [\luisD,\luisH]=-2\luisH \quad  &[\luisD,\luisC]=2\luisC 
\quad &[\luisH,\luisC]=\luisD \quad &[\luisH,\luisP]=0 \cr
 [\luisP,\luisC]= \frac 12 \luisK  \quad  &
 [\luisK,\luisH]=-2\luisP\quad  &[\luisK,\luisC]=0
  \quad  &
\end{array}
\label{bg}
\ee
that is, they come into the non-deformed Schr\"odinger Lie
algebra. Hence the deformation parameter $\sigma$  appears 
only within the coproduct, which now reads
\bea
&&\Delta(\luisM)=1\otimes \luisM + \luisM\otimes 1\cr
&&\Delta(\luisP)=1\otimes \luisP + \luisP\otimes 1 + \sigma
\luisP\otimes \luisP\cr
&&\Delta(\luisH)=1\otimes \luisH +
\luisH\otimes (1 + \sigma \luisP)^2\cr
&&\Delta(\luisK)=1\otimes
\luisK + \luisK\otimes  \frac{1}{1 + \sigma \luisP} 
-2\sigma
\luisD'\otimes \frac{\luisM}{1 + \sigma \luisP}\cr
&&\Delta(\luisD)=1\otimes \luisD + \luisD\otimes  
\frac{1}{1 +
\sigma \luisP} - 
\frac{1}2 \luisM\otimes \frac{\sigma \luisP}{1 +
\sigma \luisP}\cr
&&\Delta(\luisC)=1\otimes \luisC + \luisC\otimes
\frac{1}{(1 + \sigma \luisP)^2} -\frac \sigma 2 
\luisD'\otimes
\frac{1}{1 + \sigma \luisP}\,\luisK \cr
&&\qquad\qquad+
\frac{\sigma^2}{2}\luisD'(\luisD'-1)\otimes 
\frac{\luisM}{(1 +
\sigma
\luisP)^2}
\eea
where $ \luisD'=\luisD+\frac 12 \luisM$. Note that 
although the new generator $\luisP$ is non-primitive, its
coproduct satisifes
\be
\Delta\left((1+ {\sigma \luisP})^a\right)= (1+ {\sigma
\luisP})^a\otimes (1+ {\sigma \luisP})^a
\label{bbg}
\ee
since the old generator $P$ fulfils  $\Delta(e^{a zP})= 
e^{a zP}\otimes e^{a zP}$ for any real number $a$.

On the other hand, the mapping (\ref{bf}) transforms the 
differential-difference
realization  (\ref{bc}) into
\bea
&&
\luisP=\Delta_x \qquad \qquad \luisH=\partial_t
\qquad \qquad \luisM= m \cr
&& \luisK=2 t \Delta_x  +  m (2 x + \sigma)
T_x^{-1} \qquad \luisD=2  t  \partial_t  
+ x \Delta_x  T_x^{-1}  -
\frac 12 T_x^{-1} +1 \cr
&& \luisC= t^2 \partial_t   + t x \Delta_x
T_x^{-1}
 +\frac 12 m x^2   T_x^{-2}  
 + t \left(  1- \frac 12 T_x^{-1} \right)
-\frac {m \sigma^2}{8}  T_x^{-2}   .
\label{bh}
\eea
We emphasize that these are just the symmetry operators of 
the equation (\ref{za}) obtained by Floreanini {\it et
al.\/}${}^{9}$,  provided the
continuum time limit
$\tau\rightarrow 0$ is performed and $m=\frac 12$.   This
connection is consistent with the fact that (\ref{bf})
transforms the Casimir of the Galilei subalgebra 
$E_z^{(s)}$ (\ref{cass}) into the non-deformed one
\be
E=\luisP^2 -2 \luisM\luisH .
\label{casnd}
\ee
 Therefore, the discretized SE
\be
(\Delta_x^2- 2 m
\partial_t)\phi(x,t)=0
\label{bk}
\ee
is obtained as the realization (\ref{bh}) of $E$,  and the
operators (\ref{bh}) are symmetries of this equation 
satisfying
\be
[E,X]=0\quad
X\in\{\luisK,\luisH,\luisP,\luisM\}\qquad [E,\luisD]= 2  
E\qquad
 [E,\luisC]= 2 t  E .
\label{bl}
\ee

Thus, we  have shown  that the space 
differential-difference SE under study$\,{}^{9}$ has  
$U_\sigma^{(s)}({\cal S})$ as its quantum Hopf symmetry
algebra. The deformation parameter $\sigma=-z$ is
interpreted as the lattice step in the $x$ coordinate,
meanwhile the time $t$ remains as a continuum variable. We
also remark that, by using (\ref{bh}), the solutions of
(\ref{bk}) have been obtained$\,{}^{9}$ for $m=\frac 12$.

%%%%%%%%%%%%%%%%%%%%%%%%%%%%%%%%%%%%%%%%%%%%%%%

\section*{III. A DISCRETE TIME SCHR\"ODINGER EQUATION}

A similar procedure can be applied to the quantum 
Schr\"odinger algebra
$U_z^{(t)}({\cal S})$  coming from the non-standard 
classical
$r$-matrix
\be
r^{(t)}= 2z H\wedge (D +\frac 12 M)\equiv 2z H\wedge D'.
\label{bll}
\ee
The coproduct and commutation rules of 
$U_z^{(t)}({\cal S})$ are  given by$\,{}^{20}$
\be
\begin{array}{l}
\Delta(M)=1\otimes M + M
\otimes 1\\  [2pt]
\Delta(H)=1\otimes H + H \otimes 1 \\  [2pt]
\Delta(P)=1\otimes P + P\otimes e^{-2zH}\\  [2pt]
\Delta(K)=1\otimes K + K\otimes e^{2zH}-2z D'\otimes
e^{4zH}P\\  [2pt]
\Delta(D)=1\otimes D + D\otimes e^{4zH} +  \frac 12 M\otimes
(e^{4zH}-1)\\  [2pt]
\Delta(C)=1\otimes C + C\otimes  e^{4zH}-z   D'\otimes
e^{4zH}  M \end{array}
\label{ca}
\ee
\be
\begin{array}{l}
[D,P]=-P \qquad   [D,K]=K\qquad [K,P]=M   \qquad
[M,\,\cdot\,]=0\\  [2pt]
[D,H]= \frac 1{2z}
({1-e^{4zH}}) \qquad   
[D,C]=2 C +2z  (D')^2\qquad [H,P]=0\\ [2pt] 
[H,C]=  D' -\frac 12  M  e^{4zH}\qquad
[K,C]=z(D'K+KD') \\  [2pt]
[P,C]=-K-z  (D'P+PD')  \qquad   [K,H]= e^{4zH}P.
\end{array}
\label{cb}
\ee
A realization of (\ref{cb})   reads$\,{}^{20}$
\bea
&&H=\partial_t \qquad \qquad
P=\partial_x\qquad\qquad M=m \nonumber\\  [2pt]
&&K=- (t + 4z) e^{4 z \partial_t}
\partial_x  - m x  \qquad
  D=2 (t+ 4z) \frac{e^{4 z
\partial_t}-1}{4z} + x\partial_x + \frac 12\nonumber\\[2pt]
&&C=(t^2 - 4z b t)
\frac{e^{4 z \partial_t}-1}{4z}  + t x \partial_x 
+ \frac 12 t  + \frac
12  m x^2 - 4 z (b +1) e^{4 z \partial_t}\nonumber\\[4pt]  
&&\qquad - z x^2
\partial_{x}^2  -2 z (b+1) x \partial_x -z(b+1/2)^2 ,
\label{cc}
\eea
where $b= \frac m2 -2$. Now, a time  discretization of the 
SE is obtained by considering the deformed Casimir of the
Galilei subalgebra 
\be
E_z^{(t)} =P^2- 2 M \left( \frac{1-e^{-4zH}}{4z}\right)
\label{cast}
\ee 
written in terms of the realization (\ref{cc}):
\be
  \left(\partial_x^2 - 2 m \frac {1-e^{-4z\partial_t}}{4z}
\right)
\phi(x,t)=0 .
\label{cd}
\ee
Under the realization (\ref{cc})  the generators of
$U_z^{(t)}({\cal S})$ are symmetry operators of this 
equation as they satisfy
\bea
&&[E_z^{(t)},X]=0\quad X\in\{K,H,P,M\}\qquad [E_z^{(t)},D]=  2
E_z^{(t)} \cr
&&[E_z^{(t)},C]=2\left\{ t + z (1 - m - 2  x
\partial_x)\right\} E_z^{(t)} .
\label{ce}
\eea

The relationship with the time discretization  of the  SE
analysed by Floreanini {\it et al.\/}$\,{}^{9}$ is provided
by the non-linear map defined by
\bea
&&
\luisH=\frac   {1-e^{-4zH}}{4z} \qquad \luisP= P 
\qquad \luisM=M\cr
&&\luisD=D+ 2(1-e^{4zH}) \qquad \luisK=-2K -8z P e^{4zH}\cr
&& \luisC=C+ z (D')^2 - 4 z D e^{4zH} \qquad \tau=- 4 z .
\label{cf}
\eea
If we apply (\ref{cf}) to (\ref{cb}) we find again the 
classical commutation rules of the Schr\"odinger algebra
(\ref{bg}) while the coproduct is now given by
\bea
&&\Delta(\luisM)=1\otimes \luisM +
\luisM\otimes 1\cr
&&\Delta(\luisH)=1\otimes \luisH + \luisH\otimes 1 +
\tau \luisH\otimes \luisH\cr
&&\Delta(\luisP)=1\otimes \luisP +
\luisP\otimes (1 + \tau \luisH)^{1/2}\cr
&&\Delta(\luisK)=1\otimes
\luisK + \luisK\otimes  \frac{1}{(1 + \tau \luisH)^{1/2}} - 
\tau
\luisD'\otimes \frac{\luisP}{1 + \tau \luisH}\cr
&&\Delta(\luisD)=1\otimes \luisD
 + \luisD\otimes  \frac{1}{1 + \tau \luisH}
 - \frac{1}2 \luisM\otimes \frac{\tau \luisH}{1 + \tau
\luisH}\cr &&\Delta(\luisC)= 1\otimes \luisC + 
\luisC\otimes 
\frac{1}{1 + \tau
\luisH}   - \frac {\tau}{2} \luisD'\otimes 
  \frac{1}{1 + \tau \luisH}
\luisD \cr
&&\qquad\qquad + \frac{\tau}4 \luisD'(\luisD' - 2 )\otimes
\frac{\tau \luisH}{(1 + \tau \luisH)^2} .
\eea
The generator
$\luisH$ fulfils a similar property to (\ref{bbg}), that is,
\be
\Delta\left((1+ {\tau \luisH})^a\right)= (1+ {\tau
\luisH})^a\otimes (1+ {\tau \luisH})^a .
\ee

On the other hand, in this new basis  the
realization (\ref{cc}) turns out to be
\bea
&& \luisH=\Delta_t \qquad
\luisP=\partial_x  \qquad \luisM= m \cr
&& \luisK=2 t \partial_x
T_t^{-1}  + 2 m x \qquad \luisD=2  t  \Delta_t   T_t^{-1} + 
x
\partial_x +  \frac 12  \cr
&& \luisC= t^2 \Delta_t   T_t^{-2}    + t x
\partial_x T_t^{-1}
 +\frac 12 m x^2   + t \left( T_t^{-2} 
 - \frac 12 T_t^{-1} \right)  ,
\label{cg}
\eea
which is again the limit $\sigma\rightarrow 0$  of the 
symmetry operators of Floreanini {\it et al.\/}$\,{}^{9}$
when $m=\frac 12$. The corresponding discretized  SE is
provided by the (non-deformed) Casimir $E$ of the Galilei
subalgebra (\ref{casnd}) written through (\ref{cg}):
\be (\partial^2_x- 2 m
\Delta_t)\phi(x,t)=0 .
\label{ci}
\ee
The new operators (\ref{cg}) are
symmetries of this equation satisfying
\be
[E,X]=0\quad
X\in\{\luisK,\luisH,\luisP,\luisM\}\qquad 
[E,\luisD]=  2  E \qquad
[E,\luisC]= 2 t T_t^{-1} E .
\label{cj}
\ee
Thus we have
obtained a  discrete version of the SE on a uniform 
lattice such that its symmetry algebra is the
quantum  Hopf algebra $U_\tau^{(t)}({\cal S})$,  and the 
time  lattice step is related to the deformation parameter
in the form $\tau=-4 z$.

\subsection*{A. A NONLINEAR MAP FOR THE JORDANIAN 
DEFORMATION OF $sl(2,\R)$}

We emphasize that the four generators
$\{\luisD,\luisH,\luisC,\luisM\}$ close a quantum extended
$sl(2,\R)$ algebra.  Hence, as a byproduct, if we take
$\luisM=0$ we shall get the map
\bea
&&
\luisH=\frac   {1-e^{-4zH}}{4z} \cr
&&\luisD=D+ 2(1-e^{4zH}) \cr
&& \luisC=C+ z D^2 - 4 z D e^{4zH} \qquad \tau=- 4 z .
\eea
Now, if we denote
\be
J_3=-\luisD\qquad J_+=\luisH \qquad  J_-=-\luisC,
\ee
and we apply this map onto the well known non-standard
quantum $sl(2,\R)$ algebra$\,{}^{11-16}$
given from (\ref{ca}) and (\ref{cb})
\be
\begin{array}{l}
\Delta(H)=1\otimes H + H \otimes 1 \\  [2pt]
\Delta(D)=1\otimes D + D\otimes e^{4zH} \\  [2pt]
\Delta(C)=1\otimes C + C\otimes  e^{4zH}
\end{array}
\label{caz}
\ee
\be
[D,H]= \frac 1{2z}
({1-e^{4zH}}) \qquad   
[D,C]=2 C +2z  (D)^2 \qquad
[H,C]=  D,
\label{cbz}
\ee
we obtain a new
expression for this quantum algebra with deformed 
non-cocommutative coproduct and ``classical" commutation
rules:
\bea
&&[J_3,J_+]=2J_+\qquad
[J_3,J_-]=-2J_-\qquad [J_+,J_-]=J_3\cr
&&\Delta(J_+)=1\otimes J_+ + J_+\otimes 1 +
\tau J_+\otimes J_+\cr
&& \Delta(J_3)=1\otimes J_3
 + J_3 \otimes  \frac{1}{1 + \tau J_+} \cr 
&&\Delta(J_-)= 1\otimes J_- + J_-\otimes 
\frac{1}{1 + \tau
J_+}   + \frac {\tau}{2} J_3 \otimes 
  \frac{1}{1 + \tau J_+}
J_3 \cr
&&\qquad\qquad - \frac{\tau}4 J_3(J_3 + 2 )\otimes
\frac{\tau J_+}{(1 + \tau J_+)^2} ,
\eea
where the underlying classical $r$-matrix is 
$r=\frac{\tau}2 J_+\wedge J_3$. This result is worthy to be
compared with the previous literature on nonlinear maps for
the non-standard quantum $sl(2,\R)$ algebra of
Abdessalam {\it et al.}$\,{}^{21,22}$  since, in general,
the transformed coproduct has a very complicated form (in
this respect see the aforementioned references and 
also the work of Aizawa,${}^{23}$ where the corresponding
map is used to construct the representation theory of this
quantum algebra).

%%%%%%%%%%%%%%%%%%%%%%%%%%%%%%%%%%%%%%%%
\section*{IV. ON A POSSIBLE SPACE-TIME DISCRETIZATION}

It would be certainly interesting to analyse whether there
exists a (at least two parameter)   quantum Schr\"odinger 
algebra giving rise to the symmetries of the full 
space-time uniform discretization provided by (\ref{za}).
In general, the search for quantum deformations of a given
Lie algebra can be guided by the study of its Lie bialgebra
structures, since they are just the first order
deformations of the coproduct.${}^{24}$
Therefore,  in our case,  one could try to find the most
general Schr\"odinger bialgebra
$({\cal S},\delta)$  such that it includes among its 
components  the two Lie bialgebras linked to the previous 
deformations; hence
$({\cal S},\delta)$ should depend  on $\sigma$, $\tau$ and
perhaps on some additional parameters $\lambda_i$. Under 
such conditions and by imposing $\delta$ to  fulfil the
cocycle condition and the dual map
$\delta^\ast$ to define a Lie bracket, we  obtain a unique
cocommutator family depending only on three parameters:
$\tau=- 4 z_1 $, $\sigma=-z_2$ and $\lambda$.
  Furthermore, it can be shown
that this Schr\"odinger bialgebra is a coboundary  one
generated by the classical
$r$-matrix:
\bea
&&r=2 z_1 H\wedge D' + z_2 P\wedge D' -\ados H\wedge C
-\frac{\ados^2}{8 z_1} C\wedge D 
- \frac{2z_1 z_2}{\ados}P\wedge H\cr
&&\qquad +\frac{\ados z_2}{8 z_1}\left(K\wedge D
 -\frac{\ados}{2 z_1}
K\wedge C -3P\wedge C \right)+\left(\frac{\ados}{4} +
\frac{3z_2^2}{16 z_1}
\right) K\wedge P \cr
&&\qquad +\left( \frac{3\ados z_2}{16 z_1} +
\frac{z_2^3}{32 z_1^2}\right)
K\wedge M - \left(\frac{\ados^2}{16 z_1} 
+ \frac{\ados z_2^2}{32 z_1^2}
\right) M\wedge C \cr
&&\qquad +  \left(\frac{\ados}{4} +\frac{z_2^2}{16 z_1} 
\right) M\wedge D.
\label{sa}
\eea
We omit the explicit expressions for the  cocommutators as 
they
  can be obtained by means of $\delta(X)=[1\otimes X +
X\otimes 1,r]$.

We remark that the first two terms of (\ref{sa}) correspond 
to the time and space classical $r$-matrices  $r^{(t)}$ and
$r^{(s)}$, but  both deformations do not   fulfil a simple
``superposition principle" and  extra contributions have to  be
added. The full quantization of such Lie bialgebra seems to 
be a difficult task and we shall not address this problem
here. The ``additional" deformation parameter $\ados$
allows us to distinguish between non-standard (when $\ados
= -\frac{z_2^2}{4 z_1}$) and standard solutions
(otherwise). On the other hand, although the limits $z_2\to
0$ and $\ados\to 0$ lead to the  time classical $r$-matrix,
unfortunately the limit  $z_1\to 0$ leads to divergences in
(\ref{sa}).

%%%%%%%%%%%%%%%%%%%%%%%%%%%%%%%%%%%%%%%%
\section*{V. CONCLUDING REMARKS}

We wish to point out that the existence of a Hopf structure 
for the symmetries of a given equation associated to an
elementary system allows us to write equations of composed
systems keeping the same symmetry algebra.${}^{25,26}$
In order to use this property for the
two cases here discussed, we see
that only the last commutator in either (\ref{bl}) or 
(\ref{cj}) involving the conformal generator  $\luisC$ is
not algebraic, but depends explicitly on the chosen
representation (the same happens at the continuum level).
Consequently the composed systems characterized by the
equation
$\Delta(E) \phi = 0$ will have, by construction,
$\Delta(\luisH),\Delta(\luisP), \Delta(\luisK),
\Delta(\luisD)$, and $\Delta(\luisM)$ as symmetry 
operators  (moreover they close a Hopf subalgebra).
However, in general  this will not be the case for
$\Delta(\luisC)$, and a further study on the behaviour of
this operator is needed in order to construct coupled
equations with full Schr\"odinger algebra symmetry.

Finally, we wish to point out that the applicability of 
the constructive approach presented here is not limited to
the cases analyzed before; in fact it can be directly
extended to other quantum algebras by means of their
corresponding difference realizations.  Work on this line
is nowadays in progress.

%%%%%%%%%%%%%%%%% ACKNOWLEDGMENTS %%%%%%%%%%%

\noindent {\section*{ACKNOWLEDGEMENTS}}

This work has been partially supported by DGES 
(Projects PB94-1115 and PB95-0719)  from the Ministerio de
Educaci\'on y Cultura  (Spain), and also by the Junta de 
Castilla y Le\'on (Projects CO1/396,
CO2/197, and CO2/297).

\vskip2cm

\noindent
${}^{1}$ 
J. Lukierski, H. Ruegg and A. Nowicky, Phys.
Lett. B {\bf 293}, 344 (1992).

\noindent
${}^{2}$ 
L. Frappat and A. Sciarrino, Phys. Lett. B {\bf
347}, 1 (1995).

\noindent
${}^{3}$ 
A. Ballesteros, F.J. Herranz, M.A. del Olmo 
and M. Santander, Phys. Lett. B {\bf 351}, 137 (1995). 

\noindent
${}^{4}$
G.R.W. Quispel, H.W. Capel and R. Sahadevan, Phys. Lett. A  
{\bf 170},  379  (1992).

\noindent
${}^{5}$
D. Levi, L. Vinet and P. Winternitz, J. Phys. A  
{\bf 30},  633  (1997).

\noindent
${}^{6}$
R. Floreanini   and L. Vinet,   Lett. Math. Phys. 
{\bf  32}, 37   (1994).

\noindent
${}^{7}$
R. Floreanini and L. Vinet,   J. Math. Phys. 
{\bf 36}, 3134 (1995).

\noindent
${}^{8}$
V.K. Dobrev, H.D. Doebner  and C. Mrugalla,   
 J. Phys. A  {\bf 29}, 5909 (1996).

\noindent
${}^{9}$
  R. Floreanini,  J. Negro, L.M. Nieto  and L. Vinet,  
Lett. Math. Phys.  {\bf 36},  351 (1996).

\noindent
${}^{10}$
J.  Negro  and L.M. Nieto,   J. Phys.  A 
{\bf 29}, 1107 (1996).

\noindent
${}^{11}$ 
S. Majid, Class. Quantum Grav. {\bf 5}, (1988) 1587.

\noindent
${}^{12}$
Demidov E E, Manin Yu I, Mukhin E E and
Zhdanovich D V 1990
{\it Progr. Theor. Phys. Suppl.} {\bf 102} 203

\noindent
${}^{13}$ 
C. Ohn,  Lett. Math. Phys. {\bf 25}, 85 (1992).

\noindent
${}^{14}$
A.A. Vladimirov, Mod. Phys. Lett. A {\bf 8},  2573 (1993).

\noindent
${}^{15}$ 
A. Shariati, A. Aghamohammadi and M. Khorrami,
 Mod. Phys. Lett. A {\bf 11}, 187 (1996).

\noindent
${}^{16}$ 
A. Ballesteros and  F.J. Herranz, J. Phys. A {\bf 29}, 
L311 (1996).

\noindent
${}^{17}$
C.R. Hagen,   Phys. Rev. D  {\bf 5},  377 (1972).

\noindent
${}^{18}$
U. Niederer,  Helv. Phys. Acta  {\bf 45},  802  (1972).

\noindent
${}^{19}$
 A. Ballesteros,   F.J. Herranz and P. Parashar, 
J. Phys. A  {\bf 30},  8587 (1997).

\noindent
${}^{20}$
 A. Ballesteros,   F.J. Herranz and P. Parashar,  
Mod. Phys. Lett. A {\bf 13},  1241 (1998).

\noindent
${}^{21}$ 
B. Abdesselam, A. Chakrabarti and R. Chakrabarti,
 Mod. Phys. Lett. A {\bf 11}, 2883 (1996).

\noindent
${}^{22}$ 
B. Abdesselam, A. Chakrabarti and R. Chakrabarti,
 Mod. Phys. Lett. A {\bf 13}, 779 (1998).

\noindent
${}^{23}$
N. Aizawa,  J. Phys. A {\bf  30}, 5981 (1997).

\noindent
${}^{24}$
A.  Ballesteros   and  F.J. Herranz, 
J. Phys. A  {\bf 29},  4307 (1996).

\noindent
${}^{25}$
F. Bonechi, E. Celeghini, R. Giachetti, E. Sorace  and
M. Tarlini,  Phys Rev Lett {\bf  68},  3718 (1992).

\noindent
${}^{26}$
F. Bonechi, E. Celeghini, R. Giachetti, E. Sorace  and
M. Tarlini,   Phys Rev B {\bf  46}, 5727 (1992).

\end{document}